\documentclass[preprints,article,accept,pdftex,moreauthors]{my-mdpi}

\firstpage{1}
\pubvolume{6}
\issuenum{8}
\articlenumber{434}
\pubyear{2022}
\copyrightyear{2022}
\externaleditor{{Academic Editor:} Gani Stamov}
\datereceived{26 May 2022}
\daterevised{27 July 2022} 
\dateaccepted{7 August 2022}
\datepublished{10 August 2022}
\hreflink{https://doi.org/ 10.3390/fractalfract6080434} 
\doinum{10.3390/fractalfract6080434}

\pdfoutput=1

\Title{Study of a Fractional Creep Problem with Multiple Delays in Terms of Boltzmann's Superposition Principle}

\TitleCitation{Study of a Fractional Creep Problem with Multiple Delays in Terms of Boltzmann's Superposition Principle}

\Author{Amar Chidouh $^{1,\dagger}$\orcidA{},
Rahima Atmania $^{2,\dagger}$\orcidB{}
and Delfim F. M. Torres $^{3,}$*$^{,\dagger}$\orcidC{}}

\AuthorNames{Amar Chidouh, Rahima Atmania and Delfim F. M. Torres}

\AuthorCitation{Chidouh, A.; Atmania, R.; Torres, D.F.M.}

\address{$^{1}$ \quad Department of Mathematics, Chadli Bendjedid University, Eltarf 36000, Algeria; m2ma.chidouh@gmail.com\\
$^{2}$ \quad Department of Mathematics, Faculty of Sciences, Badji Mokhtar University, P.O. Box 12, \linebreak 
Annaba 23000, Algeria; atmanira@yahoo.fr\\
$^{3}$ \quad Center for Research and Development in Mathematics and Applications (CIDMA), \linebreak 
Department of Mathematics, University of Aveiro, 3810-193 Aveiro, Portugal}

\corres{Correspondence: delfim@ua.pt}

\firstnote{These authors contributed equally to this work.}

\abstract{We study a class of nonlinear fractional differential equations 
with multiple delays, which is represented by the Voigt creep fractional model of viscoelasticity. 
We discuss two Voigt models, the first being linear and the second being nonlinear. 
The linear Voigt model give us the physical interpretation and is associated 
with important results since the creep function characterizes the viscoelastic 
behavior of stress and strain. For the nonlinear model of Voigt,
our theoretical study and analysis provides existence and stability, 
where time delays are expressed in terms of Boltzmann's superposition principle.
By means of the Banach contraction principle, we prove existence of a unique solution 
and investigate its continuous dependence upon the initial data 
as well as Ulam stability. The results are illustrated with an example.}



\begin{document}


\section{Introduction}

In recent years, fractional calculus has attracted the attention
of many researchers~\mbox{\cite{MR4412846,MR4428365,105,106,107,108,109}}.		
Such a research axis is of great importance in many fields, especially
when dealing with memory or hereditary properties, such as 
in viscoelastic phenomena: for example, stress--strain 
in polymeric materials \cite{MR4409941,103,104}.

To describe the behavior of materials and show their viscoelastic properties, 
one uses rheological models, which are of Voigt or Maxwell type or a combination 
of these basic models; see \cite{MR4350570,30,33,101,102} and the references therein. 
On the other hand, there are many materials that are difficult to describe in rheological 
models that contain finite elements of elastic and viscous components. For that reason,
it is common to resort to the use of fractional models, which give us the opportunity 
to use few elements and, at the same time, give us an accurate description.
In order to understand the problem at hand, and as an initial idea, 
let us recall here the simplest form of the equation that governs the creep phenomenon:
\begin{equation}
\label{1111}
{\eta}x^{\prime }(t)+E x(t)=\varphi (t),
\quad x(t_{0})=x_{0},
\end{equation}
where $\eta $ is the viscosity coefficient and $E$ is the modulus of the elasticity.
In (\ref{1111}), $t_{0}$ should be chosen in a way that
for $t<t_{0}$, the material is at rest, without stress and strain.
For a given stress history $\varphi$, the strain $x$ is expressed by
\begin{equation}
\label{re}
x(t)=\int\limits_{t_{0}}^{t}e^{-\frac{t-s}{\tau }}\varphi (s)ds,
\quad \tau = \frac{\eta}{E},
\end{equation}
where the constant $\tau $ is called the retardation time. The creep function
associated with problem (\ref{1111}) is then given by
\begin{equation}
\label{re1}
k(t-t_{0})=\frac{1}{E}\left( 1-\exp\left(-\frac{t-t_{0}}{\tau}\right)\right),
\quad t\geq t_{0},
\end{equation}
which is, in general, a completely monotonic function; see, e.g., \cite{30,101}. 
This means that the viscoelastic function $k(t)$ must satisfy the inequalities
\begin{equation}
\label{creep c}
(-1)^{n}\frac{d^{n}}{dt^{n}}k(t)\geq 0,
\quad n=1,2,\ldots,
\end{equation}
assuring its monotonicity. It is noticeable from Formulas (\ref{re}) and (\ref{re1}) 
that for a given stress, the response is not instantaneous and the strain takes time.
This is well apparent when dealing with viscoelastic phenomena. As fractional operators 
are in general integral operators with a singular kernel, involving a time delay makes 
those operators the best for modeling such rheological phenomena: see \cite{m7,[m],reo}
and references therein.

Here, we study a more general creep problem, based on the Voigt model, that satisfies 
Boltzmann's superposition as expressed in the second terms of the following 
nonlinear fractional differential equation:
\begin{equation}
\label{a}
^{C}D_{0+}^{\alpha }x\left( t\right) +\lambda x\left( t\right)
=\underset{j=1}{\overset{n}{\sum }}b_{j}\left( t\right) g_{j}\left(
x\left( t-\tau_{j}\right) \right),
\quad t\in \left[ 0,T\right],
\end{equation}
where $^{c}D_{0}^{\alpha }$ denotes the Caputo derivative of order $0<\alpha <1$;
$\lambda $ is a real positive constant; $b_{j}$, $g_{j}$
are given functions; and $\tau _{j}$ represents the time delays. We take
\begin{equation}
\label{b}
x\left( t\right) =\psi \left( t\right),
\quad t\in \left[ -v,0\right],
\end{equation}
where $0<v=\underset{j=1,\ldots,n}{\max }\left\{ \tau _{j}\right\} \leq T$
and $x\left( 0\right) =\psi \left( 0\right) =0$.

A large number of mathematicians obtained several results
on the existence, uniqueness, and stability for equations with delays 
of the type described from \eqref{1111} to \eqref{creep c}
\cite{2,5,8,9,21,MR4128368,MR3244467,MR4414720}. Motivated by 
these papers, we obtain here sufficient conditions for the stability 
of our fractional delayed differential problem (\ref{a}) and (\ref{b}). 
Before that, we prove the existence and uniqueness of the solution, which is crucial 
to give a physical meaning to our problem. This contrasts with previous works, 
which do not combine a comprehensive mathematical analysis with physical interpretation.

The manuscript is structured as follows. In Section~\ref{sec:2}, we give some fundamental
results, which will be used throughout the paper. In particular, we obtain an integral
representation of our problem and then extract the fractional creep function associated 
with the linear problem. In Section~\ref{sec:3}, we use the Banach contraction principle 
to show the existence and uniqueness of the solution and then establish its continuous 
dependence with the initial data. In Section~\ref{sec:4}, we study the Ulam stability, 
including an illustrative example. We end with conclusions in Section~\ref{sec:5}.


\section{Fundamental Results on the Linear Problem}
\label{sec:2}

We begin by recalling the definition of the Mittag--Leffler function, which is an
important tool in fractional calculus and will have an important role in our paper.

\begin{Definition}
\label{prop1}
The two-parameter Mittag--Leffler function is defined by the
series expansion
\begin{equation}
\label{1.6}
E_{\alpha ,\beta }(z)=\overset{\infty }{\sum\limits_{n=0}}\dfrac{z^{n}}{
\Gamma (\alpha n+\beta )},
\quad z,\alpha , \beta \in\mathbb{C}
\text{ \ with }\textrm{Re}\alpha >0,
\end{equation}
where $\Gamma(\cdot)$ is Euler's gamma function.
In particular, $E_{\alpha }(z) = E_{\alpha ,1}(z)$.
\end{Definition}

Schneider in \cite{schneider} proved that the generalized Mittag--Leffler function
$E_{\alpha ,\beta }(-t)\ $with $t\geq 0$ is completely monotonic if and only
if $0<\alpha \leq 1$ and $\beta \geq \alpha $. In other words,
\begin{equation}
\label{res1}
(-1)^{n}\frac{d^{n}}{dt^{n}}E_{\alpha ,\beta }(-t)\geq 0,
\end{equation}
for all $n=0,1,2,\ldots$ Note that this property generalizes (\ref{creep c}).

From (\ref{res1}), we can see that for $0<\alpha \leq 1$ and $\beta \geq
\alpha $,
\begin{equation}
\label{res2}
E_{\alpha ,\beta }(-t)\leq \frac{1}{\Gamma (\beta )},
\quad t\geq 0,
\end{equation}
and the above estimation (\ref{res2}) will enable us to establish our results.

In the sequel, we use the notation $e_{\alpha }^{\lambda t}$ for the
$\alpha$-exponential function:
\begin{equation}
\label{f1}
e_{\alpha }^{\lambda t}:=t^{\alpha -1}E_{\alpha ,\alpha }(\lambda t^{\alpha}).
\end{equation}

\begin{Lemma}[See \cite{32}]
Let $\alpha ,\lambda \in \mathbb{C}$ and $\textrm{Re}(\alpha )>0$. Then,
\begin{equation}
\label{cc}
\int_{0}^{t}e_{\alpha }^{\lambda t}dt
=t^{\alpha }E_{\alpha ,\alpha+1}(\lambda t^{\alpha }).
\end{equation}
\end{Lemma}

Now, we consider our linear fractional problem as the fractional Voigt model
\begin{equation}
\label{ww}
\left\{
\begin{array}{l}
^{C}D_{0+}^{\alpha }x(t)+\lambda x(t)=\varphi (t),
\quad  t> 0,\quad \lambda >0,\\
x(0)=0.
\end{array}
\right.
\end{equation}

Viscoelastic phenomena can be expressed with integral equations or
differential ones. While differential equations are related to rheological models,
which provide a more direct physical interpretation of the viscoelastic behavior,
the integral equations are more general and suitable for theoretical work.
Indeed, using the Laplace transform, we can convert our linear
problem (\ref{ww}) into a Volterra integral equation as follows:
\begin{equation}
\label{h1}
x(t)=\int_{0}^{t}e_{\alpha }^{-\lambda (t-s)}\varphi (s)ds.
\end{equation}

\begin{Theorem}
\label{z}
Suppose that the stress ${\varphi}$ of the fractional initial value problem (\ref{ww})
is a continuous function on $[0,T]$. Then, (\ref{h1}) is the strain and the continuous
solution of (\ref{ww}). Moreover, (\ref{ww}) and (\ref{h1}) are equivalent in $C[0,T]$.
\end{Theorem}

\begin{proof}
Since $\varphi$ is a continuous function on $[0,T]$, then, 
from ([Theorem~3.24] of \cite{kilbass}), problem (\ref{ww}) is equivalent 
in $C[0,1]$ to the following Volterra integral equation of second~kind:
\begin{equation}
\label{rig}
x(t)=\frac{1}{\Gamma (\alpha )}\int\limits_{0}^{t}(t-s)^{\alpha -1}\varphi (s)ds
-\frac{\lambda}{\Gamma(\alpha )}\int\limits_{0}^{t}(t-s)^{\alpha -1}x(s)ds.
\end{equation}

Now, we apply the successive approximation method to solve (\ref{rig}).
Let
\begin{equation}
\label{fd}
x_{0}(t)=I^{\alpha }\varphi=\frac{1}{\Gamma (\alpha)}
\int\limits_{0}^{t}(t-s)^{\alpha -1}\varphi (s)ds.
\end{equation}

Substituting (\ref{fd}) into (\ref{rig}), we obtain that
\begin{eqnarray*}
x_{1} &=&I^{\alpha }\varphi -{\lambda}I^{\alpha }x_{0} \\
&=&I^{\alpha }\varphi -{\lambda}I^{2\alpha }\varphi,
\end{eqnarray*}
and
\begin{eqnarray*}
x_{2} &=&I^{\alpha }\varphi -{\lambda}I^{\alpha }x_{1}\\
&=&I^{\alpha }\varphi -{\lambda}I^{2\alpha }\varphi
+ {\lambda}^2I^{3\alpha }\varphi.
\end{eqnarray*}

Continuing this process, we obtain
\begin{eqnarray*}
x_{m}(t) &=&\frac{1}{\Gamma (k\alpha
+\alpha )}\int\limits_{0}^{t}(t-s)^{k\alpha +\alpha-1}
\sum\limits_{k=0}^{m}\left( {-\lambda}\right) ^{k}\varphi(s)ds\\
&=&\int\limits_{0}^{t}(t-s)^{\alpha
-1}\sum\limits_{k=0}^{m}\frac{(t-s)^{k\alpha }}{\Gamma (k\alpha +\alpha)}
\left({-\lambda}\right) ^{k}\varphi (s)ds,
\end{eqnarray*}
and, passing to the limit, as $m\rightarrow \infty$,
\begin{eqnarray*}
x(t) &=&\int\limits_{0}^{t}(t-s)^{\alpha
-1}\sum\limits_{k=0}^{\infty }\frac{(t-s)^{k\alpha }}{\Gamma (k\alpha
+\alpha )}\left( -{\lambda}\right) ^{k}\varphi (s)ds \\
&=&\int\limits_{0}^{t}(t-s)^{\alpha -1}E_{\alpha
,\alpha }\left( -{\lambda}(t-s)^{\alpha }\right)\varphi(s)ds\\
&=&\int_{0}^{t}e_{\alpha }^{-\lambda (t-s)}\varphi (s)ds.
\end{eqnarray*}

The proof is complete.
\end{proof}

\begin{Remark}
We draw attention to the fact that the solution to Equation (\ref{1111})
is described by the exponential function while the solution to
(\ref{ww}) is expressed in terms of the $\alpha$-exponential function,
which is a generalization of (\ref{re}) and plays an important role here.
\end{Remark}

Now, we look to the creep function $k_{\alpha }(t)$
and give the following result.

\begin{Theorem}
\label{pro}
The creep function associated with the fractional problem
(\ref{ww}) is given by
\begin{equation}
\label{flu1}
k_{\alpha }(t)=t^{\alpha }E_{\alpha ,\alpha +1}(-{\lambda }t^{\alpha }).
\end{equation}
\end{Theorem}

\begin{proof}
We use the stress--strain relation
\begin{equation}
\label{alt}
x(t)=\int_{0}^{t}k_{\alpha }(t-s)\varphi ^{\prime }(s),
\end{equation}
which is known as the Boltzmann--Volterra equation. Integrating by parts
(\ref{h1}) and using~(\ref{cc}), we obtain
\begin{eqnarray*}
x(t) &=&\int_{0}^{t}e_{\alpha }^{-\lambda (t-s)}\varphi (s)ds \\
&=&t^{\alpha }E_{\alpha ,\alpha +1}(-{\lambda }t^{\alpha })\varphi (0)
+\int_{0}^{t}(t-s)^{\alpha }E_{\alpha ,\alpha +1}({-\lambda }
(t-s)^{\alpha })\varphi ^{\prime }(s)ds,
\end{eqnarray*}
from which we can write that
\begin{equation}
\label{flu2}
k_{\alpha }(t)=t^{\alpha }
E_{\alpha ,\alpha +1}(-{\lambda }t^{\alpha }).
\end{equation}

The proof is complete.
\end{proof}
Taking $\lambda=\frac{E}{\eta }$ and evaluating 
the creep function \eqref{flu2}, we obtain
\begin{eqnarray*}
k_{\alpha }(t) &=&\frac{1}{\eta}t^{\alpha }E_{\alpha ,\alpha +1}\left( -\frac{1}{\tau }
t^{\alpha }\right) \\
&=&\frac{1}{\eta}t^{\alpha }\sum\limits_{n=0}^{\infty }\left( -\frac{1}{\tau }\right) ^{n}
\frac{t^{\alpha n}}{\Gamma (\alpha n+\alpha +1)} \\
&=&-\frac{\tau}{\eta} \sum\limits_{n=1}^{\infty }\left( -\frac{1}{\tau }\right) ^{n}
\frac{t^{\alpha n}}{\Gamma (\alpha n+1)} \\
&=&-\frac{1}{E} \left( E_{\alpha }\left( -\frac{1}{\tau }t^{\alpha }\right)
-1\right) \text{,}
\end{eqnarray*}
that is,
\begin{equation}
\label{fluage f}
k_{\alpha }(t)=\frac{1}{E} \left( 1-E_{\alpha }\left(
-\frac{1}{\tau }t^{\alpha}\right) \right).
\end{equation}

Note that for $\alpha =1$, one obtains from \eqref{fluage f} that
\begin{eqnarray*}
k_{1}(t) &=&\frac{1}{E} \left( 1
-\exp \left( -\left(\frac{t}{\tau }\right)\right) \right)\\
&=&k(t-t_{0})\text{ for }t_{0}=0\text{,}
\end{eqnarray*}
which is the creep function (\ref{re1}) associated with problem (\ref{1111}).


\section{Existence and Uniqueness Results}
\label{sec:3}

In this section, we consider the nonlinear problem (\ref{a}) and (\ref{b})  
as a generalization of the classical creep problem. Here, we will rely on
Boltzmann's superposition principle \cite{33}:
\begin{equation*}
\varphi (t)=\underset{j=1}{\overset{n}{\sum }}b_{j}\left( t\right)
g_{j}\left( x\left( t-\tau _{j}\right) \right),
\end{equation*}
which gives the stress as a sum of the delayed response to the strain.

We introduce the following assumptions:

\begin{Hypothesis}
\textbf{\emph{(H1):}} 
\emph{$\ b_{j}:\left[ 0,T\right] \rightarrow \mathbb{R}$
\emph{are continuous functions with} $B_{j}=\underset{t
\in \left[ 0,T\right] }{\sup}\left\vert b_{j}\left( t\right) \right\vert$.}
\end{Hypothesis}

\begin{Hypothesis}
\textbf{\emph{(H2):}}
\emph{$g_{j}:{C\left[ -v,T\right] }\rightarrow {C\left[ -v,T\right] }$
\emph{are Lipschitz functions, i.e., there exists} $l_{j}>0$ \emph{such that}
\begin{equation}
\label{L}
\left\Vert g_{j}(x)-g_{j}(y)\right\Vert _{C\left[ -v,T\right] }\leq
l_{j}\left\Vert x-y\right\Vert _{C\left[ -v,T\right] },
\quad  x,y\in C\left[ -v,T\right],
\end{equation}
\emph{with} $g_{j}\left( 0\right) \neq 0$.}
\end{Hypothesis}

\begin{Theorem}
\label{the1}
Assume that the assumptions (H1) and (H2) hold. If
\begin{equation}
\label{CL}
T^{\alpha }\underset{j=1}{\overset{n}{\sum }}B_{j}l_{j}<\Gamma (\alpha +1),
\end{equation}
then problem (\ref{a}) and (\ref{b})  has a unique solution
$x\in C\left(\left[ -v,T\right], \mathbb{R}\right)$.
\end{Theorem}

\begin{proof}
The solution of problem (\ref{a}) and (\ref{b})  satisfies, for
$t\in \left[ 0,T\right]$, the integral equation
\begin{equation}
\label{p}
x\left( t\right) =\int_{0}^{t}e_{\alpha }^{-\lambda (t-s)}
\underset{j=1}{\overset{n}{\sum }}b_{j}\left( s\right) g_{j}\left(
x\left( s-\tau_{j}\right) \right) ds,
\end{equation}
and
\begin{equation}
\label{p'}
x\left( t\right) =\psi \left( t\right) \ \text{\textit{for}}\mathit{\ }
t\in \left[ -v,0\right] \ \text{ with}\ x\left( 0\right) =\psi \left( 0\right) =0.
\end{equation}

Let $\psi $ be a continuous function and denote
\begin{equation*}
X=\left\{ x\in C\left( \left[ -v,T\right], \mathbb{R}\right) :
\text{\textit{\ }}x\left\vert _{\left[ -v,0\right] }=\psi \right.\right\},
\end{equation*}
which is a Banach space endowed with the $\sup$--norm
\begin{equation*}
\left\Vert \gamma \right\Vert _{C\left[ -v,T\right] }=\underset{t\in \left[
-v,T\right] }{\sup }\left\vert \gamma \left( t\right) \right\vert ,
\end{equation*}
and where $x\left\vert _{\left[ -v,0\right] }\right.$ is the restriction of the
function $x$ on $\left[ -v,0\right] $. Taking into account Theorem~\ref{z}
and the assumptions (H1)--(H2), we define the operator $P$ $:X\rightarrow X$ by
\begin{equation*}
(Px)(t)=\left\{
\begin{array}{c}
\psi \left( t\right) ,\mathit{\ \ \ \ \ \ \ \ \ \ \ \ \ \ \ \ \ \ \ \ \ \ \
\ \ \ \ \ \ \ \ \ \ \ \ \ \ \ \ \ \ \ \ \ \ \ }t\in \left[ -v,0\right] , \\
0,\ \ \ \ \ \ \ \ \ \ \ \ \ \ \ \ \ \ \ \ \ \ \ \ \ \ \ \ \ \ \ \ \ \ \ \ \
\ \ \ \ \ \ \ \ \ \ \ \ \ \ \ \ \ \ \ \ \ t=0, \\
\displaystyle \int_{0}^{t}e_{\alpha }^{-\lambda (t-s)}\underset{j=1}{\overset{n}{\sum }}
b_{j}\left( s\right) g_{j}\left( x\left( s-\tau _{j}\right) \right) ds,\ \ \
\ \ \ \ t\in \left[ 0,T\right].
\end{array}
\right.
\end{equation*}

Now, using the contraction principle mapping of Banach, we investigate the
existence and uniqueness of the fixed point of the operator $P$ in $X.$ In
fact, for $x,y\in X$, one has
\begin{eqnarray*}
\left\vert Px\left( t\right) -Py\left( t\right) \right\vert  &\leq
&\int_{0}^{t}e_{\alpha }^{-\lambda (t-s)}\underset{j=1}{\overset{n}{\sum }}
\left\vert b_{j}\left( s\right) \right\vert l_{j}\left\vert x(s-\tau
_{j})-y(s-\tau _{j})\right\vert ds \\
&\leq &\underset{j=1}{\overset{n}{\sum }}\int_{0}^{t}e_{\alpha }^{-\lambda
(t-s)}\left\vert b_{j}\left( s\right) \right\vert l_{j}\left\vert x(s-\tau
_{j})-y(s-\tau _{j})\right\vert ds \\
&\leq &\underset{j=1}{\overset{n}{\sum }}\int_{-\tau _{j}}^{0}e_{\alpha
}^{-\lambda (t-z-\tau _{j})}\left\vert b_{j}\left( s\right) \right\vert
l_{j}\left\vert x(z)-y(z)\right\vert dz \\
&&+\underset{j=1}{\overset{n}{\sum }}\int_{0}^{t-\tau _{j}}e_{\alpha
}^{-\lambda (t-z-\tau _{j})}\left\vert b_{j}\left( s\right) \right\vert
l_{j}\left\vert x(z)-y(z)\right\vert dz.
\end{eqnarray*}
As $x(t)=y(t)=\psi \left( t\right) $ for $\mathit{\ }t\in \left[ -v,0\right]$,
we have
\begin{equation*}
\left\vert Px\left( t\right) -Py\left( t\right) \right\vert
\leq \underset{j=1}{\overset{n}{\sum }}\int_{0}^{t-\tau _{j}}
e_{\alpha }^{-\lambda(t-z-\tau _{j})}\left\vert b_{j}\left( s\right)
\right\vert l_{j}\left\vert x(z)-y(z)\right\vert dz,
\end{equation*}
and, using (\ref{res2}) and (\ref{cc}), we obtain
\begin{equation*}
\left\Vert Px-Py\right\Vert _{X}\leq \ \left( \frac{T^{\alpha }}{\Gamma
(\alpha +1)}\underset{j=1}{\overset{n}{\sum }}B_{j}l_{j}\right) \left\Vert
x-y\right\Vert _{X},
\end{equation*}
which shows by (\ref{CL}) that $P$ is a contraction in $X$. Thus, from the Banach
fixed point theorem, we conclude that $P$ has a unique fixed point in $C\left[
-v,T\right]$, which is the unique solution of (\ref{a}) and (\ref{b})  in
$C\left( \left[-v,T\right],\mathbb{R}\right)$.
\end{proof}

As a consequence of Theorem~\ref{the1},
we prove the continuous dependence of the solution
with respect to the initial data of the problem.

\begin{Corollary}
\label{the}
Under the conditions of Theorem \ref{the1}, the unique solution
of (\ref{a}) and (\ref{b}) depends continuously
on function $\psi \left( t\right)$.
\end{Corollary}

\begin{proof}
Let $x_{1}$and $x_{2}$ be solutions of Equations (\ref{a}) and (\ref{b}) 
corresponding to the initial data $\psi _{1}\left( t\right) $ 
and $\psi _{2}\left( t\right)$, respectively. Then,
\begin{equation*}
\left\vert x_{1}\left( t\right) -x_{2}\left( t\right) \right\vert
=\int_{0}^{t}e_{\alpha }^{-\lambda (t-s)}\underset{j=1}{\overset{n}{\sum }}
\left\vert b_{j}\left( s\right) \right\vert l_{j}\left\vert x_{1}\left(
s-\tau _{j}\right) -x_{2}\left( s-\tau _{j}\right) \right\vert ds.
\end{equation*}

Putting $\left( s-\tau _{j}\right) =z$, we obtain
\begin{equation*}
\begin{split}
\left\vert x_{1}\left( t\right) -x_{2}\left( t\right) \right\vert
&\leq \frac{1}{\Gamma (\alpha )}\overset{n}{\sum_{j=1}}B_{j}l_{j}\int
\limits_{-\tau _{j}}^{t-\tau _{j}}\left( t-\tau _{j}-z\right) ^{\alpha
-1}\left\vert x_{1}\left( z\right) -x_{2}\left( z\right) \right\vert dz \\
&\leq \frac{1}{\Gamma (\alpha )}\overset{n}{\sum_{j=1}}B_{j}l_{j}\int
\limits_{-\tau _{j}}^{0}\left( t-\tau _{j}-z\right) ^{\alpha -1}\left\vert
x_{1}\left( z\right) -x_{2}\left( z\right) \right\vert dz \\
&\quad +\frac{1}{\Gamma (\alpha )}\overset{n}{\sum_{j=1}}B_{j}l_{j}\int
\limits_{0}^{t-\tau _{j}}\left( t-\tau _{j}-z\right) ^{\alpha -1}\left\vert
x_{1}\left( z\right) -x_{2}\left( z\right) \right\vert dz \\
&\leq \frac{1}{\Gamma (\alpha )}\overset{n}{\sum_{j=1}}B_{j}l_{j}\sup_{z
\in \lbrack -v,0]}\left\vert x_{1}\left( z\right) -x_{2}\left( z\right)
\right\vert \int\limits_{-\tau _{j}}^{0}\left( t-\tau _{j}-z\right)
^{\alpha -1}dz \\
&\quad +\frac{1}{\Gamma (\alpha )}\overset{n}{\sum_{j=1}}B_{j}l_{j}\sup_{z\in
\lbrack 0,T]}\left\vert x_{1}\left( z\right) -x_{2}\left( z\right)
\right\vert \int\limits_{0}^{t-\tau _{j}}\left( t-\tau _{j}-z\right)
^{\alpha -1}dz \\
&\leq \frac{1}{\Gamma(\alpha +1)}\overset{n}{\sum_{j=1}}
\left(t^{\alpha }-\left( t-\tau _{j}\right) ^{\alpha }\right)
B_{j}l_{j}\left\Vert \psi _{1}-\psi_{2}\right\Vert _{C\left[ -v,0\right]}\\
&\quad +\frac{1}{\Gamma (\alpha +1)}\overset{n}{\sum_{j=1}}
\left( t-\tau _{j}\right) ^{\alpha }
B_{j}l_{j}\left\Vert x_{1}-x_{2}\right\Vert _{C\left[ 0,T\right]}.
\end{split}
\end{equation*}

Consequently, we obtain that
\begin{equation*}
\left\Vert x_{1}-x_{2}\right\Vert _{C\left[ 0,T\right] }
\leq \frac{T^{\alpha }}{\Gamma (\alpha +1)}\overset{n}{\sum_{j=1}}B_{j}l_{j}\left\Vert
x_{1}-x_{2}\right\Vert _{C\left[ 0,T\right] }
+\frac{T^{\alpha }}{\Gamma (\alpha +1)}\overset{n}{\sum_{j=1}}
B_{j}l_{j}\left\Vert \psi _{1}-\psi _{2}\right\Vert _{C\left[ -v,0\right] }.
\end{equation*}

In view of (\ref{CL}), one has
\begin{equation}
\label{ab}
\left\Vert x_{1}-x_{2}\right\Vert _{C\left[ 0,T\right] }\leq \frac{T^{\alpha
}\overset{n}{\sum_{j=1}}B_{j}l_{j}}{\Gamma (\alpha +1)-T^{\alpha }\overset{n}{
\sum_{j=1}}B_{j}l_{j}}\left\Vert \psi _{1}-\psi _{2}\right\Vert _{C\left[
-v,0\right] }.
\end{equation}
This implies the continuous dependence of $x$ on the initial data $\psi$.
\end{proof}


\section{Ulam--Hyers Stability}
\label{sec:4}

One of the main qualitative properties of solutions of differential
equations is stability, which is studied by many methods \cite{MR4232864,MR3980195}.
Recently, Ulam-type stabilities have attracted more and more attention
\cite{MR4413902,MR4407931}. The classical concept of Ulam
stability was posed by Ulam in 1940 and later
obtained for functional equations by Hyers in 1941 \cite{10}.
Hyers' result was extended by replacing functional equations with
differential equations, and this approach guarantees the existence of an
$\varepsilon$-solution, which is quite useful in many applications where
finding the exact solution is impossible. Now, let us give the definition of
Ulam--Hyers stability in the fractional setting. For more details,
we refer the reader to \cite{2,8,10} and references therein.

\begin{Definition}
The fractional differential equation
\begin{equation}
\label{eq}
\Psi (t,y,f,^{c}D^{\alpha _{1}},\ldots,^{c}D^{\alpha _{n}})=0
\end{equation}
is Hyers--Ulam stable if, for a given $\varepsilon >0$ and a function
$y_{\varepsilon }$ such that
\begin{equation*}
\left\vert \Psi (t,y_{\varepsilon },f,^{c}D^{\alpha _{1}},
\ldots,^{c}D^{\alpha _{n}})\right\vert
\leq \varepsilon,
\end{equation*}
there exists a solution $y_{e}$ of (\ref{eq}) and a positive constant $K>0$
such that
\begin{equation*}
\left\vert y_{\varepsilon }(t)-y_{e}(t)\right\vert \leq K\varepsilon .
\end{equation*}
\end{Definition}

\begin{Theorem}
\label{hi}
If problem (\ref{a}) and (\ref{b}) has a unique solution, then it is Ulam--Hyers stable.
\end{Theorem}

\begin{proof}
Let $y_{e}$ be a unique solution in $C\left( \left[ -v,T\right],
\mathbb{R}\right)$ satisfying (\ref{p}) and (\ref{p'}) and $y_{\varepsilon }$ be a
solution of \ the following inequality:
\begin{equation}
\label{eq1}
\left\vert ^{C}D_{0+}^{\alpha }y_{\varepsilon }\left( t\right)
+\lambda y_{\varepsilon}\left( t\right)
-\underset{j=1}{\overset{n}{\sum }}b_{j}\left( t\right)
g_{j}\left( y_{\varepsilon }\left( t-\tau _{j}\right) \right)\right\vert
\leq \varepsilon,
\quad t\in \left[ 0,T\right],
\end{equation}
which means there exists a function $h(t)$ such that $\left\vert h\left(
t\right) \right\vert \leq \varepsilon $ for every $t\in \left[ 0,T\right]$.

Hence, for a given continuous function $\psi$, we have
\begin{equation*}
\left\{
\begin{array}{c}
^{C}D_{0^{+}}^{\alpha }y_{\varepsilon }\left( t\right) +\lambda
y_{\varepsilon }\left( t\right) =\underset{j=1}{\overset{n}{\sum }}
b_{j}\left( t\right) g_{j}\left( y_{\varepsilon }\left( t-\tau _{j}\right)
\right) +h\left( t\right),
\quad t\in \left[ 0,T\right],  \\
y_{\varepsilon }\left( t\right) =\psi \left( t\right),
\quad t\in \left[-v,0\right].
\end{array}
\right.
\end{equation*}

In addition, we have
\begin{equation}
\label{ffff}
\left\vert y_{\varepsilon }\left( t\right) -\int_{0}^{t}e_{\alpha
}^{-\lambda (t-s)}\underset{j=1}{\overset{n}{\sum }}b_{j}\left( t\right)
g_{j}\left( y_{\varepsilon }\left( s-\tau _{j}\right) \right)ds \right\vert
\leq \frac{T^{\alpha }}{\Gamma (\alpha +1)}\varepsilon.
\end{equation}

Keeping in mind the above inequality (\ref{ffff}), we obtain
\begin{eqnarray*}
\left\vert y_{\varepsilon }\left( t\right) -y_{e}\left( t\right) \right\vert
&\leq &\left\vert y_{\varepsilon }\left( t\right) -\int_{0}^{t}e_{\alpha
}^{-\lambda (t-s)}\underset{j=1}{\overset{n}{\sum }}b_{j}\left( t\right)
g_{j}\left( y_{e}\left( s-\tau _{j}\right) \right) ds\right\vert  \\
&\leq &\left\vert y_{\varepsilon }\left( t\right) -\int_{0}^{t}e_{\alpha
}^{-\lambda (t-s)}\underset{j=1}{\overset{n}{\sum }}b_{j}\left( t\right)
g_{j}\left( y_{\varepsilon }\left( s-\tau _{j}\right) \right) ds\right\vert
\\
&&+\left\vert \int_{0}^{t}e_{\alpha }^{-\lambda (t-s)}
\underset{j=1}{\overset{n}{\sum }}b_{j}\left( t\right)
\left\vert g_{j}\left( y_{\varepsilon}\left( s-\tau _{j}\right)
-g_{j}(y_{e}\left( s-\tau _{j}\right) \right)
\right\vert ds\right\vert  \\
&\leq &\frac{T^{\alpha }}{\Gamma (\alpha +1)}\varepsilon +\frac{1}{\Gamma
(\alpha )}\overset{n}{\sum_{j=1}}B_{j}l_{j}\int\limits_{-\tau
_{j}}^{t-\tau _{j}}(t-\tau _{j}-z)^{\alpha -1}\left\vert y_{\varepsilon
}\left( z\right) -y_{e}\left( z\right) \right\vert dz.
\end{eqnarray*}

As $y_{\varepsilon }\left( z\right) =y_{e}\left( z\right) =\psi (z)$, for
$z\in \lbrack -\tau _{j},0]$, we obtain
\begin{eqnarray*}
\left\vert y_{\varepsilon }\left( t\right) -y_{e}\left( t\right) \right\vert
&\leq &\frac{T^{\alpha }}{\Gamma (\alpha +1)}\varepsilon +\frac{1}{\Gamma
(\alpha )}\overset{n}{\sum_{j=1}}B_{j}l_{j}\int\limits_{0}^{t-\tau
_{j}}(t-\tau _{j}-z)^{\alpha -1}\left\vert y_{\varepsilon }\left( z\right)
-y_{e}\left( z\right) \right\vert dz \\
&\leq &\frac{T^{\alpha }}{\Gamma (\alpha +1)}\varepsilon +\frac{T^{\alpha
}\sum\limits_{j=1}^{n}B_{j}l_{j}}{\Gamma (\alpha +1)}\left\Vert
y_{\varepsilon }-y_{e}\right\Vert _{_{C\left[ 0,T\right] }},
\end{eqnarray*}
and hence
\begin{equation*}
\left\Vert y_{\varepsilon }-y_{e}\right\Vert _{C\left[ 0,T\right] }
<\frac{T^{\alpha }}{\Gamma (\alpha +1)-T^{\alpha }\sum\limits_{j=1}^{n}B_{j}l_{j}}
\varepsilon.
\end{equation*}

In this way, we find a constant $K=\frac{T^{\alpha }}{\Gamma (\alpha
+1)-T^{\alpha }\sum\limits_{j=1}^{n}B_{j}l_{j}}$, which is well defined
taking into account condition (\ref{CL}).
\end{proof}

We end with an example of application of our results.

\begin{Example}
Consider the following creep fractional problem:
\begin{equation}
\label{px}
\left\{
\begin{array}{c}
^{C}D_{0+}^{\frac{1}{2}}x\left( t\right) +x\left( t\right)
=\sum\limits_{j=1}^{n}t^{j}\frac{x\left( t-\frac{1}{j}\right) +j}{j+3},
\quad t\in \left[ 0,1\right] , \\
x\left( t\right) =t,
\quad t\in \left[ -1,0\right],
\end{array}
\right.
\end{equation}
where $\alpha =\frac{1}{2},\ \lambda =1,\ b_{j}(t)=t^{j}$ and
$g_{j}(x(t-\tau _{j}))=\frac{x\left( t-\frac{1}{j}\right) +j}{j+3}$
with $\tau_{j}=\frac{1}{j}$ for $j=\overline{1,3}$.

First, we know from Theorem~\ref{pro} that the creep function $k_{\alpha }(t)$
associated with the linear problem of (\ref{px}) is given by
\begin{equation*}
k_{\frac{1}{2}}(t)=\left( 1-E_{\frac{1}{2}}\left( -t^{\frac{1}{2}}\right)
\right).
\end{equation*}

Hence,
\begin{equation*}
k_{\frac{1}{2}}(t)=1-e^{t}\left( 1+\textrm{erf}\left( -t^{\frac{1}{2}}\right)
\right) ,
\end{equation*}
where
\begin{equation*}
\textrm{erf}\left( z\right) =\frac{2}{\sqrt{\pi }}\int
\limits_{0}^{z}e^{-t^{2}}dt.
\end{equation*}

We have that $g_{j}(x(t-\tau _{j}))=\frac{x\left( t-\frac{1}{j}\right) +j}{j+3}$
are Lipschitz functions with $l_{j}=\frac{1}{j+3}$. Then, for $T=1$,
$\alpha=\frac{1}{2}$, $n=3$, $l_{1}=\frac{1}{4}$, $l_{2}=\frac{1}{5}$, $l_{3}=\frac{1}{6}$
and $B_{j}=1$ for $j=\overline{1,3}$, we obtain
\begin{equation*}
T^{\alpha }\sum\limits_{j=1}^{n}B_{j}l_{j}=0.62
\leq \Gamma\left(\frac{3}{2}\right)\approx 0.87,
\end{equation*}
and condition (\ref{CL}) holds. It follows from our Theorem~\ref{the1} that
problem (\ref{px}) has a unique solution in
$C\left( \left[ -1,1\right],\mathbb{R}\right)$. Finally,
by Theorem~\ref{hi}, we conclude that problem (\ref{px}) is Ulam--Hyers stable.
\end{Example}


\section{Conclusions}
\label{sec:5}

In this paper, we have studied a class of nonlinear fractional differential equations 
with multiple delays as a Voigt model, which are expressed by Boltzmann's superposition principle. 
We first took into account the Voigt linear model, because it expresses the physical meaning, 
giving strength to our problem and explaining its connection with real-life models 
where stability analysis is one of the most important investigation topics. A new form 
to our model is then given as a Volterra integral equation involving a creep kernel 
as a generalized exponential function. Such a kind of integral equation has shown 
to be very appropriate to our theoretical work, allowing the existence and stability analysis. 
In particular, using the fixed point method and Banach's contraction mapping principle, 
we succeeded to give a sufficient condition for establishing the important result of existence. 
Then, we investigated the continuous dependence upon the initial data and Ulam's stability.


\vspace{6pt}


\authorcontributions{Conceptualization, A.C., R.A. and D.F.M.T.;
methodology, A.C., R.A. and D.F.M.T.;
validation, A.C., R.A. and D.F.M.T.;
formal analysis, A.C., R.A. and D.F.M.T.;
investigation, A.C., R.A. and D.F.M.T.;
writing---original draft preparation, A.C., R.A. and D.F.M.T.;
writing---review and editing, A.C., R.A. and D.F.M.T.
All authors have read and agreed to the published version of \mbox{the manuscript}.}

\funding{This research was partially funded by FCT, grant number UIDB/04106/2020 (CIDMA).}

\institutionalreview{Not applicable.}

\informedconsent{Not applicable.}

\dataavailability{Not applicable.}

\acknowledgments{The authors are grateful to three reviewers for several
constructive comments, suggestions and questions that helped 
them to improve their work.}

\conflictsofinterest{The authors declare no conflict of interest.
The funder had no role in the design of the study; in the collection,
analyses, or interpretation of data; in the writing of the manuscript,
or in the decision to publish the~results.}


\end{paracol}

\reftitle{References}


\end{document}